\documentclass {article}

\usepackage {amsmath}
\usepackage {amssymb}
\usepackage {color}
\usepackage {graphicx}

\usepackage {geometry}
\newcommand {\boxrel} {\ensuremath{\mathbin{\,\square\,}}}

\newcommand {\prf}      {\noindent \textit{Proof.} }
\newcommand {\prfn} [1] {\noindent \textit{Proof of #1.} }
\newcommand {\qed}      {\hfill $\Box$}
\newcommand {\qedd}     {\hfill $\Box$\bigskip}

\newtheorem {defn}        {Definition} [section]
\newtheorem {lem}  [defn] {Lemma}
\newtheorem {thm}  [defn] {Theorem}
\newtheorem {cor}  [defn] {Corollary}

\numberwithin {equation} {section}

\newenvironment{romanlist}
  {

   \begin{enumerate}
  }{
   \end{enumerate}
  }

\title  {The Cartesian product of graphs with loops}
\author {Tetiana Boiko \thanks{Tetiana Boiko, Johannes Cuno, Wilfried Imrich, and Florian Lehner are supported by the Austrian Science Fund (FWF): W1230, Doctoral Program ``Discrete Mathematics.''} \\ Technische Universit\"{a}t Graz \\ \texttt{boiko@math.tugraz.at} \and Johannes Cuno \\ Technische Universit\"{a}t Graz \\ \texttt{cuno@math.tugraz.at} \and Wilfried Imrich \\ Montanuniversit\"{a}t Leoben \\ \texttt{imrich@unileoben.ac.at} \and Florian Lehner \\ Technische Universit\"{a}t Graz \\ \texttt{f.lehner@tugraz.at} \and Christiaan E. van de Woestijne \thanks{Christiaan E. van de Woestijne is supported by the Austrian Science Fund (FWF): S9611. This project is part of the Austrian National Research Network ``Analytic Combinatorics and Probabilistic Number Theory.''} \\ Montanuniversit\"{a}t Leoben \\ \texttt{c.vandewoestijne@unileoben.ac.at}}

\begin{document}

\maketitle

\begin{abstract}
We extend the definition of the Cartesian product to graphs with loops and show that the Sabidussi--Vizing unique factorization theorem for connected finite simple graphs still holds in this context for all connected finite graphs with at least one unlooped vertex. We also prove that this factorization can be computed in $O(m)$ time, where $m$ is the number of edges of the given graph. \smallskip


 \smallskip

%

\noindent \textbf{MSC classes:} 05C70, 13A05, 20M13, 05C85.
\end{abstract}


\section{Introduction}

This paper considers finite undirected graphs that may contain loops, or, put differently, symmetric binary relations on finite sets. One may define several binary operations on such graphs; these are explored in the recently revised monograph \cite{haim-2011}. The well-known \emph{Cartesian product} of finite undirected graphs is usually defined only for \emph{simple} graphs, that is, for graphs that do not contain multiple edges between the same pair of vertices and, more importantly for us, do not contain loops. Here we extend this definition.

Before doing so, let us fix the notation. For us, a \emph{graph} $G=(V,E)$ will always be a finite undirected graph without multiple edges. The edge set $E$ is taken to be a set of ordered pairs of vertices; thus, a loop on the vertex $v\in V$ corresponds to the edge $(v,v)\in E$, and as all graphs are undirected, we have $(v,w)\in E$ if and only if $(w,v)\in E$. We will occasionally call a loop a \emph{$1$-edge} and an edge that is not a loop a \emph{$2$-edge}. Moreover, given a graph $G$, we will refer to its vertex set as $V(G)$ and to its edge set as $E(G)$.

\begin{defn}[Cartesian product] \label{DefCart}
Let $G_1,\ldots,G_k$ be graphs. The Cartesian product $G=G_1\boxrel \cdots \boxrel G_k$ is a graph with vertex set $V(G)=V(G_1) \times \cdots \times V(G_k)$, and edge set $E(G)$ defined as follows: two vertices $(v_1,\ldots,v_k)\in V(G)$ and $(w_1,\ldots,w_k)\in V(G)$ are adjacent if there exists an index $i$ such that $(v_i,w_i)\in E(G_i)$, and $v_j=w_j$ for all $j\ne i$.
\end{defn}

Note that this definition extends the classical one for simple graphs. The product graph has a loop on a vertex $(v_1,\ldots,v_k)\in V(G)$ if and only if there is a loop on at least one of the constituents $v_i\in V(G_i)$. Thus, the distribution of loops (or $1$-edges) on the product graph is independent from the distribution of the $2$-edges.

\begin{defn}[projection] \label{DefProj}
Let $G_1,\ldots,G_k$ be graphs, and $G=G_1\boxrel\cdots\boxrel G_k$. The $i$th projection $p_i:V(G)\rightarrow V(G_i)$ is given by $(v_1,\ldots,v_k)\mapsto v_i$.
\end{defn}

Using Definition~\ref{DefCart}, we preserve the property that the projections $p_i:V(G)\rightarrow V(G_i)$ are weak homomorphisms from $G$ to $G_{i}$. Recall that a \emph{weak homomorphism} between graphs $G$ and $H$ is a map $\varphi:V(G)\rightarrow V(H)$ such that, whenever $(v,w)\in E(G)$, either $(\varphi(v),\varphi(w))\in E(H)$ or $\varphi(v)=\varphi(w)$. In particular, the presence of loops in $G$ or $H$ does not impose any restriction on a weak homomorphism from $G$ to $H$.

\begin{defn}[layer] \label{DefLayer}
Let $G_1,\ldots,G_k$ be graphs, and $G=G_1\boxrel\cdots\boxrel G_k$. For every vertex $a=(a_1,\ldots,a_k)\in V(G)$, the $G_i$-layer through $a$ is the induced subgraph
\[
\begin{array}{r@{\,}c@{\,}l}
G_i^a & = & \langle \{ x\in V(G) \mid p_j(x)=a_j \text{ for } j\ne i\} \rangle \smallskip \\
& = & \langle \{ (a_1,a_2,\ldots,x_i,\ldots,a_k) \mid x_i\in V(G_i)\} \rangle.
\end{array}
\]
\end{defn}

Note that $G_i^a=G_i^b$ if and only if $p_j(a)=p_j(b)$ for each index $j\ne i$. With the usual Cartesian product, the restrictions $p_i|V(G_i^a):V(G_i^a) \rightarrow V(G_i)$ are isomorphisms between $G_i^a$ and $G_i$ \cite[Section\ 4.3]{haim-2011}. Under Definition~\ref{DefCart}, we obtain a dichotomy, as follows.

\begin{lem} \label{LemDicho}
Let $G_1,\ldots,G_k$ be graphs, and $G=G_1\boxrel\cdots\boxrel G_k$. Then, the following two conditions hold for every vertex $a=(a_1,\ldots,a_k)\in V(G)$ and every $i\in\{1,\ldots,k\}$:

\begin{enumerate}
\item[(i)] If $a_j\in V(G_j)$ is unlooped for every $j\ne i$, then $p_i|V(G_i^a):V(G_i^a) \rightarrow V(G_i)$ is an isomorphism between $G_i^a$ and $G_i$.
\item[(ii)] Otherwise, $G_i^a$ is isomorphic to $G_i$ with a loop attached to every vertex.
\end{enumerate}
\end{lem}

\prf Easy from the definitions. \qed


\section{Matrix and semiring properties}

From the definition of the Cartesian product we infer that it is commutative and distributive over the disjoint union. Moreover, the trivial graph $K_1$, that is, a vertex without edges, is a unit. As the Cartesian product is also associative, see below, the set $\Gamma_0$ of isomorphism classes of finite undirected graphs with loops is a commutative semiring.

To prove associativity we could adapt the proof of \cite[Proposition 4.1]{haim-2011} for associativity of the Cartesian product of graphs without loops, or we could modify the multiplication table method of \cite[Exercise 4.15]{haim-2011}, which was introduced for the classification of associative products. However, we follow a different path and use the fact that the adjacency matrix $A(G\boxrel H)$ of the Cartesian product of two simple graphs is the Kronecker sum of the adjacency matrices $A(G)$ and $A(H)$ of the factors, see \cite[Section 33.3]{haim-2011}.

Let us first recall that the Kronecker sum $A\oplus B$ of an $n \times n$ matrix ${A}$ by an $m \times m$ matrix $B$ is defined as $I_n\otimes B + A\otimes I_m$. Here, $I_n$ and $I_m$ denote the identity matrices of size $n$ and $m$, respectively, and $P\otimes Q$ denotes the Kronecker product. In our situation, the first factor $P=(p_{ij})$ is always an $n\times n$ matrix and the Kronecker product is defined by
\[
{P}\otimes {Q} =\begin{bmatrix} p_{11} {Q} & \cdots & p_{1n}{Q} \\ \vdots & \ddots & \vdots \\ p_{n1} {Q} & \cdots & p_{nn} {Q} \end{bmatrix}.
\]
Notice that both the Kronecker sum and the Kronecker product are associative but not commutative.

For simple graphs $G$ and $H$ we have $A(G\boxrel H) = A(G)\oplus A(H)$. For graphs with loops we find that the diagonal entries take positive integer values that are not restricted to $\{0,1\}$. If we agree on the convention that a positive diagonal entry in the adjacency matrix means a loop, whereas a $0$ means no loop, then the product given in Definition~\ref{DefCart} still corresponds to the Kronecker sum. It follows that, up to isomorphism of graphs, this product is associative.

We note in passing that the fact that the Kronecker sum is not commutative does not contradict the commutativity of the Cartesian product: $A(G)\oplus A(H)$ and $A(H)\oplus A(G)$ represent adjacency matrices of $G\boxrel H$ for different vertex numberings.

Finally, we briefly call a graph \emph{entirely looped} if every vertex has a loop. For any graph $G$, we let ${\cal N}(G)$ be $G$ with its loops removed.

\begin{lem} \label{LemLooped}
Let $G$, $H$, $H_1$, $H_2$ be graphs. Assume that $G$ is entirely looped. Then $G\boxrel H$ is entirely looped as well. Moreover, if ${\cal N}(H_1) \cong {\cal N}(H_2)$, then $G \boxrel H_1 \cong G \boxrel H_2$.
\end{lem}

\prf The first statement follows directly from Definition~\ref{DefCart}. As remarked earlier, the $2$-edges of the products $G\boxrel H_i$ do not depend on the loops of either factor. Thus
\[
\begin{array}{r@{\;}c@{\;}l} {\cal N}(G\boxrel H_1) & = & {\cal N}(G) \boxrel {\cal N}(H_1) \smallskip \\ & \cong & {\cal N}(G) \boxrel {\cal N}(H_2) \smallskip \\ & = & {\cal N}(G \boxrel H_2)\,.\end{array}
\]
Next, we insert the loops on the product; but, as every vertex of $G$ has a loop, it follows that every vertex of either product $G\boxrel H_i$ has a loop as well, and the two products are obviously isomorphic. \qedd

It follows that the subset $\Gamma_{00}$ of $\Gamma_0$ given by the isomorphism classes of entirely looped graphs constitutes an \emph{ideal} of the semiring $\Gamma_0$. It is obviously closed under the disjoint union and the Cartesian product, and, since the loop $K_1^*$ is a unit for the Cartesian product inside $\Gamma_{00}$, it is a semiring itself. The loop-removing map ${\cal N}$ constitutes an isomorphism of semirings between $\Gamma_{00}$ and the set of simple graphs $\Gamma$.


\section{Unique factorization}

One fundamental property of the Cartesian product, proved independently by Sabidussi \cite{sa-1960} and Vizing \cite{vi-1963} in the 1960s, is the unique factorization of connected simple graphs into irreducibles with respect to this product. We will extend this result to graphs with loops, where we will have to exclude the set of entirely looped graphs (Lemma~\ref{LemLooped} suggests why). Algebraically speaking, we might want to form the \emph{quotient semiring} $\Gamma_0/\Gamma_{00}$, so that also any fully looped components in disconnected graphs are annulled. However, since we will only consider connected graphs in what follows, this is not of great consequence.

\begin{defn}[irreducible] \label{DefIrred}
A nontrivial, connected graph $G$ with at least one unlooped vertex is called irreducible with respect to the Cartesian product if, for every factorization $G= H \boxrel L$, either $H$ or $L$ is trivial.
\end{defn}

Recall that a graph is called \emph{trivial} if it is a vertex without edges. Consider a nontrivial, connected graph $G$ with at least one unlooped vertex. One can easily check that, if $G$ is not irreducible, it can be expressed as Cartesian product of two factors each of which is, again, a nontrivial, connected graph with at least one unlooped vertex. Iteration of this procedure yields a representation of $G$ as a product of irreducible graphs. It is occasionally called a \emph{prime factorization}.

Another way to prove the existence of a prime factorization is the following: Any factorization of $G$ with a maximum number of nontrivial factors must be a product of irreducible graphs. If $G$ has $n$ vertices, this maximum number is at most $\log_2 (n)$.

Our main results are the following.

\begin{thm} \label{ThmUniqueF}
Every nontrivial, connected graph with at least one unlooped vertex has a representation as a product of irreducible graphs with respect to the Cartesian product. The representation is unique up to isomorphisms and the order of the factors.
\end{thm}

\begin{thm} \label{FindF}
The unique prime factorization with respect to the Cartesian product of a nontrivial, connected graph $G$ with at least one unlooped vertex can be computed in $O(m)$ time, where $m$ is the number of edges of $G$.
\end{thm}

To prove Theorem \ref{ThmUniqueF}, we follow the method of \cite[Section 6.1]{haim-2011}, for Theorem \ref{FindF} we extend the ideas of \cite{impe-2007}. First, let us define convex subgraphs and boxes.

\begin{defn}[convex and box]
A subgraph $H$ of a graph $G$ is convex in $G$ if every shortest path in $G$ that connects two vertices of $H$ is completely contained in $H$. A subgraph $H$ of a Cartesian product $G=G_1\boxrel \cdots\boxrel G_k$ is called a box or subproduct if there are subgraphs $H_i\subseteq G_i$ such that
\[
H = H_1\boxrel\cdots\boxrel H_k\,.
\]
\end{defn}

In order to determine whether a subgraph is convex or not, only the $2$-edges need to be concerned. In particular, a subgraph $H$ is convex in $G$ if and only if the subgraph ${\cal N}(H)$ is convex in ${\cal N}(G)$.

\begin{lem} \label{le:convex}
Let $H$ be a subgraph of a Cartesian product $G=G_1\boxrel\cdots\boxrel G_k$. Then the following are equivalent:

\begin{enumerate}
\item[(i)] $H$ is an induced and convex subgraph of $G$;
\item[(ii)] There are induced and convex subgraphs $H_i\subseteq G_i$ such that $H=H_1\boxrel\cdots\boxrel H_k$. In other words, $H$ is a box whose factors are induced and convex.
\end{enumerate}
\end{lem}

\prf As far as only the $2$-edges are concerned, all convex subgraphs are induced and the assertion is Lemma~6.5 of \cite{haim-2011}. This means that {$p_1(V(H))\times\cdots\times p_k(V(H)) =  V(H)$}. Now, let $H_i$ be the subgraph of $G_i$ induced by $p_i(V(H))$, where $i \in \{1,\ldots,k\}$. Then the lemma follows by the definition of the Cartesian product. \qedd

\begin{figure}
\begin{center}
\input{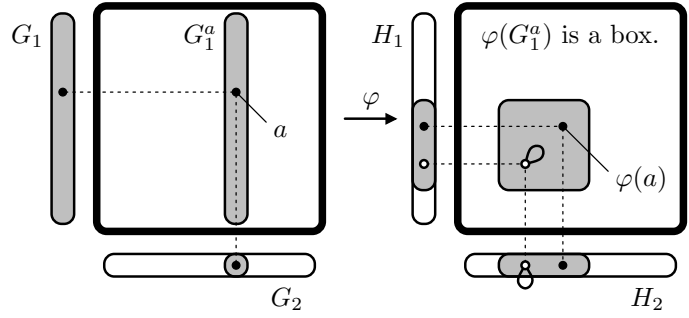}
\caption{An isomorphism between factored graphs with loops.} \label{Fig1}
\end{center}
\end{figure}

As remarked after Definition~\ref{DefIrred} every finite graph has a factorization into irreducibles. Thus we only have to show that it is unique in order to prove Theorem \ref{ThmUniqueF}. The next lemma and its corollary makes this precise; the situation is illustrated in Figure~\ref{Fig1}.


\begin{lem} \label{LemPerm}
Let $\varphi$ be an isomorphism between nontrivial, connected graphs $G$ and $H$ with at least one unlooped vertex. Assume that $G$ and $H$ are representable as products $G = G_1 \boxrel \cdots \boxrel G_k$ and $H = H_1 \boxrel \cdots \boxrel H_\ell$ of irreducible graphs. Then $k=\ell$ and, for every unlooped vertex $a\in V(G)$, there is a permutation $\pi$ of $\{1,\ldots,k\}$ such that
\[
\varphi( G_i^a ) = H_{\pi(i)}^{\varphi(a)} \text{ for every } i\in\{1,\ldots,k\}\,.
\]
\end{lem}

Formally, $\varphi$ is a bijection between the vertex sets $V(G)$ and $V(H)$. But since $\varphi$ is a homomorphism of graphs, it induces a well-defined mapping between the edge sets $E(G)$ and $E(H)$. In the above theorem, we slightly abuse notation and denote the image of the subgraph $G_i^a$, including vertices and edges, by $\varphi(G_i^a)$.
\bigskip

\prf Fix an unlooped vertex $a=(a_1,\ldots,a_k)\in V(G)$, and set $(b_1,\ldots,b_{\ell}) := \varphi(a)$. By Lemma~\ref{LemDicho} we infer that $G_i^a\cong G_i$ and $H_j^{\varphi(a)}\cong H_j$ for every $i$ and $j$. Every layer $G_i^a$ is induced and, as a consequence of Lemma~\ref{le:convex}, convex in $G$. So, its image $\varphi(G_i^a)$ is induced and convex in $H$. Again, as a consequence of  Lemma~\ref{le:convex}, $\varphi(G_i^a) = U_1 \boxrel \cdots \boxrel U_{\ell}$, where every $U_j$ is induced and convex in $H_j$. But $\varphi(G_i^a) \cong G_i^a \cong G_i$ is irreducible. Since $(b_1,\ldots,b_{\ell}) = \varphi(a) \in \varphi(G_i^a)$, we conclude that $V(U_j)=\{b_j\}$ for all indices but one, say $\pi(i)$. In other words, $\varphi(G_i^a) \subseteq H_{\pi(i)}^{\varphi(a)}$. But then
\[
G_i^a \subseteq \varphi^{-1}\left( H_{\pi(i)}^{\varphi(a)} \right)\,.
\]
Because the latter graph is induced and convex, it is a box; and because it is irreducible, it must be contained in $G_i^a$. Therefore, $\varphi(G_i^a) = H_{\pi(i)}^{\varphi(a)}$.

We claim that the map $\pi:\{1,\ldots,k\}\rightarrow \{1,\ldots,\ell\}$ is injective. If $\pi(i)=\pi(j)$, then
\[
\varphi( G_i^a) = H_{\pi(i)}^{\varphi(a)} = H_{\pi(j)}^{\varphi(a)} = \varphi(G_j^a)\,.
\]
But $\varphi$ is an isomorphism, and therefore the above equation implies $G_i^a=G_j^a$. Since every layer contains at least two vertices, we obtain $i=j$. So, $\pi$ is injective, and $k\le \ell$. Repetition of the above argument for $\varphi^{-1}$ yields $\ell\le k$. So, $k=\ell$ and $\pi$ is a permutation.
\qed

\begin{cor} \label{cor:main}
$G_i \cong H_{\pi(i)}$ for every $i\in\{1,\ldots,k\}$.
\end{cor}

\prf Since $a$ is unlooped, $G_i \cong G_i^a$ and $H_j \cong H_j^{\varphi(a)}$ for every $i$ and $j$. By Lemma \ref{LemPerm} the corollary follows. \qedd

Clearly Lemma \ref{LemPerm} and Corollary \ref{cor:main} prove the validity of Theorem \ref{ThmUniqueF}.

\subsection*{A remark about automorphisms}

In Lemma \ref{LemPerm} the permutation $\pi$ of $\{1,\ldots,k\}$ is constructed to a fixed unlooped vertex $a\in V(G)$. Actually $\pi$ is independent of the choice of $a$, and one can extend  Lemma \ref{LemPerm} to the following description of the automorphisms of $G$.

\begin{thm} \label{Auto} Suppose $\varphi$ is an automorphism of a nontrivial, connected graph $G$ with at least one unlooped vertex and prime factorization $G = G_1\boxrel \cdots \boxrel G_k$. Then there are a permutation $\pi$ of $\{1,\ldots,k\}$ and isomorphisms $\varphi_i : G_{\pi(i)} \rightarrow  G_i $ for which
\[
\varphi(x_1, \ldots , x_k) = (\varphi_1(x_{\pi(1)}), \ldots , \varphi_k(x_{\pi(k)}))\,.
\]
\end{thm}

The proof of this theorem can be led on the same lines as that of \cite[Theorem 6.10]{haim-2011}. Among other consequences this implies that the automorphism group of $G$ is isomorphic to the automorphism group of the disjoint union of the prime factors $G_1, \ldots, G_k$.


\section{Algorithms}

In this section we present two algorithms for the decomposition of a nontrivial, connected graph $G$ with at least one unlooped vertex into its prime factors. One is straightforward and has complexity $O(mn)$, where $m$ is the number of edges and $n$ the number of vertices of $G$. The other one is linear in the number of edges of $G$ and depends on the algorithm of Imrich and Peterin \cite{impe-2007} for the prime factorization of graphs without loops.

Let $G = G_1\boxrel \cdots \boxrel G_k$ be the prime factorization of a nontrivial, connected graph $G$ with at least one unlooped vertex. Then also $\mathcal{N}(G) = \mathcal{N}(G_1)\boxrel \cdots \boxrel \mathcal{N}(G_k)$. Clearly the graphs $\mathcal{N}(G_i)$, $i\in\{1,\ldots,k\}$, need not be irreducible with respect to the Cartesian product. Let $\mathcal{N}(G_i) = H_{i,1}\boxrel \cdots \boxrel H_{i, \ell(i)}$ be their prime factorizations. Thus
\[
\mathcal{N}(G) = \prod_{i=1}^{k}\prod_{j=1}^{\ell(i)}H_{i,j}
\]
is a representation of $\mathcal{N}(G)$ as a Cartesian product of irreducible graphs. Because the prime factorization is unique, it is the prime factorization of $\mathcal{N}(G)$, up to the order and isomorphisms of the factors. In other words, if $\prod_{j\in J}Z_j$ is a prime factorization of $\mathcal{N}(G)$, then there is a partition $J =J_1\cup\cdots\cup J_k$ such that $\mathcal{N}(G_i) = \prod_{j\in J_i}Z_j$. Our task is to find this partition. We begin with a straightforward approach and prove the following lemma.

\begin{lem}
Let $G$ be a nontrivial, connected graph with at least one unlooped vertex. Then its prime factorization can be found in $O(mn)$ time.
\end{lem}

\prf If $G$ has $n$ vertices, then this is also true for $\mathcal{N}(G)$, and so the number of factors of $\mathcal{N}(G)$, say $r$, is at most $\log_2(n)$. This also bounds the size of $J$ and implies that the number $s$ of subsets of $J$ is at most $2^{\log_2(n)},$ i.\,e.~$s\leq n$. Notice that the factors of {$\mathcal{N}(G)$} can be found in $O(m)$ time by \cite{impe-2007}.

Let  $J_1, J_2, \ldots, J_{s}$ be all subsets of $J$, ordered in such a way that $|J_i| \leq |J_j|$ whenever $1 \leq i \leq j \leq s$. For every $i\in\{1,\ldots,s\}$ set  $Y_i := \prod_{j\in J_i}Z_j$ and $Y_i^\ast := \prod_{j\in J\setminus J_i}Z_j$. Let $\langle Y_i^a\rangle_G$ denote the subgraph of $G$ {induced} by the layer $Y_i^a$ of $Y_i$ through $a$, and define $\langle (Y_i^\ast)^a\rangle_G$ analogously. If the partition $J_i \cup (J\setminus J_i)$ of $J$ leads to a factorization of $G$,  then $\langle Y_i^a\rangle_G$ is isomorphic to a factor of $G$.

We begin the algorithm by scanning the $J_i$ in the given order. For every $J_i$ and every vertex $v\in V(G)$ we consider the projections $p_{Y_i}(v)$ and $p_{Y_i^\ast}(v)$ into $\langle Y_i^a\rangle_G$ and $\langle (Y_i^\ast)^a\rangle_G$. If $v = (v_1, \ldots, v_r)$, then $p_{Y_i}(v) = (w_1, \ldots, w_r)$, where  $w_j = v_j$ if $j \in J_i$, and $w_j = a_j$ otherwise. Notice that $p_{Y_i}(v)$ is the vertex of shortest distance from $v$ in $\langle Y_i^a\rangle_G$. The other projection $p_{Y_i^\ast}(v)$  is defined analogously. Again, $p_{Y_i^\ast}(v)$ is the vertex of shortest distance from $v$ in $\langle (Y_i^\ast)^a\rangle_G$. Clearly $G = \langle Y_i^a\rangle_G \boxrel \langle (Y_i^\ast)^a\rangle_G$ if and only if for every vertex $v\in V(G)$ the following two conditions are satisfied:

\begin{enumerate}
\item  If $v$ is unlooped, then both $p_{Y_i}(v)$ and $p_{Y_i^\ast}(v)$ are unlooped.
\item  If $v$ has a loop then at least one of the vertices $p_{Y_i}(v)$, $p_{Y_i^\ast}(v)$ has a loop.
\end{enumerate}

\noindent The time necessary to compute $p_{Y_i}(v)$ and $p_{Y_i^\ast}(v)$ for a given $v$ is proportional to $r$. As one can check in constant time whether $p_{Y_i}(v)$ or $p_{Y_i^\ast}(v)$ has a loop, one can check in $O(nr)$ time whether $G = \langle Y_i^a\rangle_G \boxrel \langle (Y_i^\ast)^a\rangle_G$.

Notice that $r$ is the number of factors of $\mathcal{N}(G)$, which is also bounded by the minimum degree $\delta$ of $\mathcal{N}(G)$. This is easily seen, since every vertex meets every layer and, in a connected graph,  is incident with at least one edge of that layer. Hence the number of factors cannot exceed the degree of any vertex, and $nr \leq n\delta \leq m$.

For a given  $J_i$ one can thus check in $O(m)$ time whether {$\langle Y_i^a\rangle_G$} is a factor of $G$. If it is, and if $J_i$ is minimal with respect to inclusion, then it clearly is an irreducible factor. Hence, this is true for the first factor that we encounter, because of having ordered the $J_i$ by size. We now continue the scan, omitting the $J_j$ that are not disjoint from $J_i$, to find the next factor. Clearly it will also be irreducible. We continue until we have found all irreducible factors. Since there are no more than $n$ subsets of $J$, we can find them in $O(nm)$ time.\qedd

In order to reduce the complexity to $O(m)$, we need some more preparation. So let $a$ be an unlooped vertex of $G$ and  $L_i$ be the levels of a BFS-ordering of the vertices of $G$ with respect to the root $a$. That is, $L_i$ consists of all vertices of distance $i$ from $a$. Furthermore, we enumerate the vertices of $G$ by giving them so-called BFS-numbers that satisfy BFS$(v) > $ BFS$(u)$ if the distance from $a$ to $v$ is larger than the one from $a$ to $u$.

It is important to observe that the projection $p_{Y_i}(v)$ is a vertex of $\langle Y_i^a\rangle_G$ and always closer to $a$ than $v$, unless $v$ already is a vertex of  $\langle Y_i^a\rangle_G$, because then $p_{Y_i}(v) = v$. \bigskip

\prfn{Theorem \ref{FindF}} Let $\prod_{j\in J}Z_j$  be a prime factorization of $\mathcal{N}(G)$. We begin with the trivial partition of $J$ and wish to check, whether it already leads to a factorization of $G$. We scan the vertices $v$ of $G$ in BFS-order and, given $v$, check the validity of Conditions

\begin{romanlist}
\item \label{c1} If $v$ is unlooped, then  all {$p_{Y_i}(v)$} are unlooped.
\item \label{c2} If $v$ has a loop, then at least one of the projections $p_{Y_i}(v)$ has a loop.
\end{romanlist}

\noindent If one of these conditions is not satisfied, then the partition of $J$ is obviously inconsistent with the loop structure. In either case we have too many factors and have to make the partition of $J$ coarser. Before we go on, notice that in $L_1$ these conditions are trivially satisfied for any partition of $J$, because all projections $p_{Y_i}(v)$ are $a$, except one, which is $v$.

Suppose we arrive at a vertex $v$ where one of the conditions (\ref{c1}) or (\ref{c2}) is violated for the first time. Assume first that Condition (\ref{c1}) is violated, that is, $v$ is unlooped, but $p_{Y_i}(v)$ has a loop for an index $i$. In the end, all projections have to be unlooped. We must combine the set $J_i$ with one or more other sets of the partition. Using the fact that we proceed in BFS-order, it is easy to see that we have to make $v$ a unit layer vertex, that is, we combine all those sets $J_j$ for which $p_{Y_j}(v) \neq a$. Assume now that Condition (\ref{c2}) is violated, that is, $v$ has a loop, but no $p_{Y_i}(v)$ does. In the end, at least one of the projections has to have a loop. As above, the only way to achieve this is to make $v$ a {unit layer vertex}, that is, we combine all factors $J_j$ for which $p_{Y_j}(v) \neq a$.

In both cases we arrive at a coarser partition of $J$ than the one we started out with. By associativity of the Cartesian product with loops, we need not recheck the vertices we have already considered and continue in BFS-order.

Notice that this process yields a factorization, because both (\ref{c1}) and (\ref{c2}) are satisfied. For every finer partition of $J$ {one of these conditions is violated}, hence the factorization is the unique prime factorization we are looking for.

Considering the computational cost of these operations, we observe that all projections that we need for the $n$ vertices can be computed,  in $O(n|J|)$ time. Since we can check in constant time whether a vertex has a loop or not, the checks for  conditions (\ref{c1}) and (\ref{c2}) can also be done in $O(n|J|)$ time. As $|J| \leq \delta$, we have $O(n|J|) = O(n\delta) = O(m)$.

Finally, {recomputing the partition} needs at most $O(|J|)$ time, and this has to be done only at most $|J|$ times, so the cost is $O(\delta^2)$.\qed


\section{Remarks}

In \cite{im-1967} it was shown that connected set systems, or hypergraphs, as they are called now, also have unique prime factorizations with respect to the Cartesian product if one-element sets, or loops in our terminology, are excluded. Our result also extends to hypergraphs with loops: Connected hypergraphs have unique prime factorization with respect to the Cartesian product, if there is a least one vertex without a loop. Furthermore, the same arguments yield unique prime factorization for connected infinite graphs or hypergraphs with respect to the weak Cartesian product; compare \cite{im-1971}.


\end{document}